\documentclass[hidelinks,onefignum,onetabnum]{siamart220329}

\usepackage{amsfonts,upquote}
\usepackage{graphicx,verbatim}
\headers{AAA rational approximation on a continuum}{Driscoll, Nakatsukasa, and Trefethen}

\title{AAA rational approximation on a continuum
\thanks{Submitted to the editors DATE.}}

\author{Toby Driscoll\thanks{Dept.\ of Mathematical Sciences, U. of Delaware,
Newark, DE 19716, USA}
  (\email{driscoll@udel.edu}).
\and Yuji Nakatsukasa\thanks{Mathematical Institute, U. of Oxford, Oxford OX4 4DY, UK}
  (\email{nakatsukasa@maths.ox.ac.uk}).
\and Lloyd N. Trefethen\thanks{Mathematical Institute, U. of Oxford, Oxford OX4 4DY, UK}
  (\email{trefethen@maths.ox.ac.uk}).}

\def\real{{\mathbb{R}}}
\def\complex{{\mathbb{C}}}
\def\Re{\hbox{\rm Re\kern .5pt}}
\def\Im{\hbox{\rm Im\kern .5pt}}
\def\snew{s_{\hbox{\scriptsize new}}}

\begin{document}

\maketitle

\begin{abstract}
AAA rational approximation has normally been carried out on a discrete
set, typically hundreds or thousands of points in a real interval
or complex domain. Here we introduce a continuum AAA algorithm
that discretizes a domain adaptively as it goes.  This enables
fast computation of high-accuracy rational approximations on
domains such as the unit interval, the unit circle, and the
imaginary axis, even in some cases where resolution of singularities
requires exponentially clustered sample points, support points,
and poles.  Prototype MATLAB (or Octave) and Julia codes {\tt aaax}, {\tt aaaz},
and {\tt aaai} are provided for these three special domains; the
latter two are equivalent by a M\"obius transformation.
Execution is very fast since the matrices whose SVD\kern .5pt s
are computed have only three times as many rows as columns.  The codes
include a AAA-Lawson option for improvement of a AAA approximant to
minimax, so long as the accuracy is well above machine precision.
The result returned is pole-free in the approximation domain.
\end{abstract}

\begin{keywords}
AAA algorithm, rational approximation, minimax, MATLAB, Julia
\end{keywords}

\begin{MSCcodes}
41A20, 65D15
\end{MSCcodes}

\section{Introduction}
The AAA algorithm, introduced in 2018~\cite{aaa}, is a numerical
method for rational approximation of a function $f$ on a real or complex domain.
Because of its speed, reliability, and domain flexibility,
AAA has become a standard method for these computations,
with a rapidly growing literature.  Applications to date
include analytic continuation~\cite{lustri,analcont},
interpolation of equispaced data~\cite{equi}, Laplace
problems with applications to magnetics~\cite{costa2,costa}, conformal mapping~\cite{cm1,cm2},
Stokes flow~\cite{xue}, simulation of turbulence~\cite{kkw},
nonlinear eigenvalue problems~\cite{gnt,lmpv,sem},
finite element linearizations~\cite{djm}, design
of preconditioners~\cite{budisa,dpp}, model order
reduction~\cite{abg,gg,johns1,johns,pradovera}, and signal
processing~\cite{dpp,hochman,RF,vre,wdt}.
For signal processing, AAA is the basis of the
{\tt rational} code in the MathWorks RF Toolbox~\cite{RF},
and there have also been generalizations to multivariate
approximation~\cite{cbg,gg,hochman,johns1,lmpv}.

In some applications, one wants a rational approximation on
a discrete set~$Z$ of $N$ points of $\real$ or $\complex$.
More often, however, one would like to approximate on a continuum
$E$ such as the unit interval $[-1,1]$, the unit circle or unit
disk, or the imaginary axis.  In such cases AAA has normally been
applied by approximating $E$ by a fixed discrete set~$Z$, chosen
in advance, with $N$ typically in the hundreds or thousands.
The algorithm computes repeated SVD\kern .5pt s (singular value
decompositions) involving Loewner matrices with about $N$ rows
and $m = 1,2,3,\dots$ columns, where $m-1$ is the degree of the
rational function being constructed.  Since $m$ is usually below
$100$ and the operation count of such an SVD grows only linearly
with $N$, the cost associated with having many more rows than
columns is usually not too great.

Nevertheless, it would be good to have a AAA algorithm that works
directly with the continuum.  We see three reasons for this.
First, it is unappealing philosophically to require the user to
specify a discrete set rather than the domain that is actually
of interest.  Second, although the cost of user discretization
may only be a constant factor (loosely speaking), it is still
a shame if this factor is large, in cases where
the Loewner matrix is very tall and skinny.  This may become
particularly important in situations where the evaluation
of $f$ is expensive.  Third is the challenge
of computing approximations with poles and zeros clustering
near the approximation domain $E$.  Exponential clustering
of poles near singularities is essential, for it is the
source of the special power of rational approximation in
such cases~\cite{newman,clustering}.  But for this to work,
exponentially clustered sample points are needed too, and
where should they be placed?  As experienced AAA users,
we have become adept at defining grids by expressions
like \verb|logspace(-14,0,1000)| (for approximation
on $[\kern .3pt 0,1]$ with a singularity at $0$) and
\verb|tanh(linspace(-16,16,1000))| (for approximation on $[-1,1]$
with singularities at both ends), but some experimentation is
always needed to get it right, and the user ingenuity has to
increase when there are more singularities.  And, of course,
some problems have singularities at unknown locations, as can
occur in model order reduction applications when there are poles
near the imaginary axis.

This paper introduces an algorithm for continuum AAA approximation
along with template \hbox{MATLAB} (or Octave)
and \hbox{Julia} codes, each about 100 lines long,
for the special cases of the unit interval, the unit circle or
disk (these two are distinct for reasons we shall discuss),
and the imaginary axis or right half-plane.  The algorithm proceeds by greedily
choosing support points in the usual AAA fashion, which then
determine a new sample grid at every step defined by having three
equispaced sample points between each pair of support points.  Just six
new sample points are added at each step.  (One could reduce
the number to four, but experiments suggest this is
less robust.)  Poles of the current approximation are computed
at every AAA step, and if ``bad'' poles appear (i.e., poles in
$E\kern 1pt$), the current approximation is never returned as the final
AAA approximant.  Following the AAA-Lawson variant algorithm
introduced in~\cite{aaalawson}, the user may specify that the AAA
iteration should be followed by a barycentric Lawson iteration
to improve the approximation to minimax form.

\begin{figure}
\begin{center}
\vskip 10pt
\includegraphics[scale=.75]{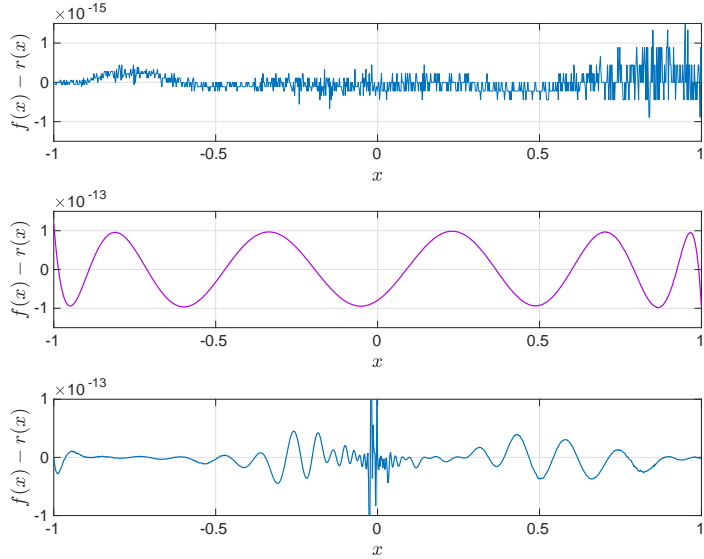}
\end{center}
\caption{\label{fig1}Illustration of\/ {\tt aaax} for rational
approximation $r\approx f$ on $[-1,1]$, with each image showing
the absolute error $f(x)-r(x)$ sampled at $1001$ equispaced points.  Top:
AAA approximation of $f(x) = e^x$ to default accuracy $10^{-13}$
(\/$1$ ms on our desktop, degree $6$).  Middle: AAA-Lawson minimax
approximation of $f(x) = e^x$ of degree $5$ (\/$5$ ms).  Here and
in other figures, AAA-Lawson error curves are distinguished by a
purple color.  Bottom: AAA approximation of $f(x) = |x|$ to requested
accuracy $10^{-13}$ actually terminates with accuracy 
$1.3\times 10^{-12}$; see also Figures $\ref{fig7}$ and
$\ref{fig8}$ (\/$600$ ms, degree $110$).} \end{figure}

\begin{figure}
\begin{center}
\vskip 10pt
\includegraphics[scale=.9]{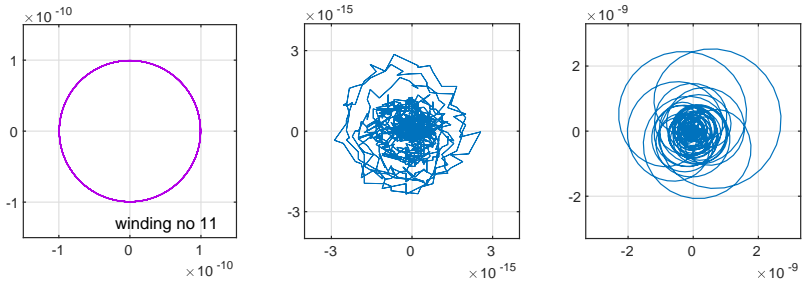}
\end{center}
\caption{\label{fig2}Illustration of\/ {\tt aaaz} for approximation
on the unit circle, with each image showing the error curve
$f(z)-r(z)$, where $z$ ranges over $1000$ equispaced points
on the circle.  Left: AAA-Lawson minimax approximation of
$f(z) = e^z$ of degree $5$ (\/$12$ ms, error curve of winding
number $11$).  Middle: AAA approximation of $\sqrt{1-z^{-2}/4}$
to default accuracy $10^{-13}$ (\/$5$ ms, degree $12$).  Right:
AAA approximation of $\sqrt{1-z}$ again with default accuracy specification
$10^{-13}$, though the code terminates at $10^{-9}$
(\/$110$ ms, degree $56$).  In this last case the
parameter {\tt mero} is set to $0$, so the approximation is
forced to be analytic in the unit disk.} \end{figure}

Figures~\ref{fig1} and~\ref{fig2} give an indication of
the behavior of the continuum AAA algorithm and the template
codes.  (These and all our illustrations are based on our MATLAB rather than
Julia codes, except in Figure~\ref{fig17}.)
The first image of Figure~\ref{fig1} shows the error in
approximation of $e^x$ on $[-1,1]$ 
to the default tolerance of $10^{-13}$.  The function
call \verb|r = aaax(@exp)| delivers a function handle \verb|r|
corresponding to a barycentric representation of a rational
approximation of $f$ of degree 6, the computation taking less than
a millisecond on a laptop (excluding plotting).  The second line
shows the error in minimax approximation of $e^x$
computed by \verb|r = aaax(@exp,5,20)|, which specifies rational
degree 5 and 20 steps of AAA-Lawson iteration (our usual number).
Note our convention of using a purple color to distinguish error
curves obtained from AAA-Lawson computations.  This computation takes
about 5 ms; the Chebfun {\tt minimax} code takes about 8 times
as long.  (All the timings in this paper are approximate.)
The third line treats the classic example~\cite{newman}
of a more difficult function with a singularity, $f(x)
= |x|$.  Here the final degree is 110, corresponding to a
much heavier computation, but still the execution is fast,
about half a second. 
The approximation returned has maximal
error about $1.3\times 10^{-12}$, which is not far from the
theoretical minimum for an approximation of this degree, about
$4\times 10^{-14}$~\cite{stahl}.\footnote{For a polynomial to
approximate $|x|$ to accuracy $1.3\times 10^{-12}$ on $[-1,1]$,
it would need degree about $215{,}000{,}000{,}000$~\cite{vc}.}
By contrast the Chebfun {\tt minimax} code~\cite{minimax} fails
above degree $82$, even though this is the most powerful Remez
algorithm implementation available, and for this it requires
four minutes to compute an approximation with the optimal error
$3.1\times 10^{-12}$.

All these computations benefit crucially from the numerical
stability of the bary\-centric representation of $r$, as is used
by both AAA approximation and {\tt minimax}.  Varga, Ruttan,
and Carpenter showed that if one works
with the quotient representation $r(x) = p(x)/q(x)$, where $p$
and $q$ are polynomials, then computing a degree 80 approximation
of $|x|$ requires 200-digit arithmetic~\cite{vrc}.  The explanation of
this instability is given in the discussion of Figure~\ref{fig5} in
the next section; see also~\cite{minimax,aaa}.

Figure~\ref{fig2} gives a similar illustration of \verb|aaaz|
for complex approximation on the unit circle.  The first image
shows degree 5 minimax approximation of $e^z$, with an error
curve $(f-r)(|z|=1)$ that looks perfectly circular.  Near-circularity of error
curves is a general phenomenon in complex minimax approximation,
and for this example the exact error curve is circular
to about one part in $10^{12}$~\cite{nearcirc,nm}.  The second
image shows AAA approximation to the default tolerance of $f(z) =
\sqrt{1-z^{-2}/4}$, a function that is analytic on the unit circle
but has branch points at $z =\pm 1/2$.  For this approximation the
flag {\tt mero}, to be discussed in Section~3,
has been set to $1$ to allow poles in the unit disk; thus we
are truly approximating on the circle, not the disk.  The third
image shows AAA approximation
of $f(z) = \sqrt{1-z}$, which has a singularity
on the unit circle.  This time {\tt mero} has been left at its
default value $0$, forcing the approximation to be analytic in
the unit disk.  The reason the accuracy 
is only $10^{-9}$ is that after this point, the
AAA iteration produced a string of approximations with
poles in the unit disk, so the code
has reverted to its best approximation found so far that is pole-free in the disk,
as will be discussed in Section~2.

The next four sections present the details of continuum AAA on
the unit interval, the unit circle or disk, the imaginary axis
or right half-plane, and other real and complex domains.
Many examples are shown along the way.  Section 6 discusses
certain challenges faced by the algorithm.
The final section summarizes various aspects of continuum AAA
approximation and its prospects.  Our current
algorithm is certainly not the last word on this subject,
but it appears to be a substantial
advance over what has been available before.

We would like to conclude this introduction with some comments on
the relationship between minimax approximations, which can be computed
by the AAA-Lawson variant, and near-minimax approximations,
computed by AAA without Lawson.

Minimax approximations are undeniably fascinating.  In real
approximation, they are characterized by equioscillatory error
curves, and in complex approximation they feature the complex
analogue shown in Figure~\ref{fig2}, error curves that are nearly
circular, sometimes so close to circular that the variation
of modulus could not be detected in floating point arithmetic.
This makes minimax approximations beautiful, memorable, and easily
recognized at a glance.

Partly for this reason, there has been a strong emphasis
in the field of approximation theory on best as opposed to near-best
approximations.  This emphasis has been amplified by the fact
that approximation theory, although it has always had
its numerical practitioners, has been dominated by theorists.

We argue that for applications, minimax approximation should
not normally be the starting point.  The basic reason is one of
practicality. Nobody knows how to compute minimax approximations
as quickly and reliably as near-minimax ones.  In the case of AAA
approximation, adding a Lawson phase typically roughly doubles
the computation time for the benefit of typically gaining just
about a digit of accuracy, and with a greater risk of failure.
(The size of the gain can be judged by noting the gaps between the
green circles and the blue dots above them in Figures~\ref{fig4},
\ref{fig5}, \ref{fig9}, and~\ref{fig11} below.)  In particular,
AAA-Lawson will almost always fail if one is dealing with
approximations of accuracy close to machine precision.  Thus in
practice, as in the first two panes of Figure~\ref{fig1}, one
is often forced to choose between a non-optimal approximation
of accuracy close to machine precision and an ``optimal''
approximation of lower accuracy!  (The latter has lower degree.)

To put it another way, even if we keep well away from machine
precision, the practical difference between minimax and
near-minimax approximations is not that the former achieves
slightly better accuracy, but that it achieves a prescribed
accuracy with a slightly lower degree.  This benefit is a
modest one, easily outweighed by a reduction in speed and
robustness.

\section{Continuum AAA on \hbox{\boldmath $[-1,1]$}}
The precise specification of our algorithm can be found in
the MATLAB code of the Appendix, particularly the
$32$ lines labeled ``Main AAA loop.''  In this section we
discuss the essential features, assuming the reader already has
some familiarity with AAA~\cite{aaa}.  

The highest level description of the algorithm is as follows.
At each step we have a row vector $S= [s_1,\dots,s_m]$ of $m$ support points in
$[-1,1]$ and a column vector $X = [x_1,\dots,x_N]^T$ of $N = 3(m-1)$ sample points
in $[-1,1]\backslash S$,
three equispaced sample points between each pair of support points.
The sample points are constructed from $S$ by a function $X =
\hbox{\tt XS}(S)$.

\vskip 8pt

{\leftskip=28pt\parindent=0pt\em\obeylines
\baselineskip 13pt
$S = [-1,1]$
{\bf for} $m = 2, 3,\dots$ {\bf until} convergence
~~~$(1)$ $X := \hbox{\tt XS}(S)$
~~~$(2)$ Use SVD to compute barycentric weights for next approximation $r\approx f$
~~~$(3)$ $S := S \cup \{\hbox{a sample point $x_i\in X$ where $|f(x_i)-r(x_i)|$ is maximal\/}\}$
{\bf end}
\par}

\vskip 8pt

Step (2) is the linearized least-squares
computation introduced in~\cite{aaa},
involving an $N\times m$ Loewner matrix $A$ with entries
\begin{equation} 
a_{ij} = {f(x_i)-f(s_j) \over x_i - s_j}.
\end{equation} 
If the $m$-vector $w$ is a minimal singular vector of $A$, then
\begin{equation}
\|Aw\|_2^{} = \hbox{minimum}, \quad \|w\|_2^{}=1,
\label{Aw}
\end{equation}
or equivalently,
\begin{equation}
\|fd - n\|_2^{} = \hbox{minimum}, \quad \|w\|_2^{}=1,
\label{fdn}
\end{equation}
where $n$ and $d$ are the numerator and denominator of the barycentric
quotient
\begin{equation}
r(x) = \sum_{j=1}^m {w_j f(s_j)\over x-s_j} \left/ \sum_{j=1}^m {w_j\over x-s_j} \right. .
\label{theapprox}
\end{equation}
In both (\ref{Aw}) and (\ref{fdn}), $\|\cdot\|_2$ is the discrete
$2$-norm over $X$.

Step (3) amounts to greedy selection of the next support point,
the standard nonlinear step of AAA approximation.\footnote{We are
struck by an analogy.  One of the jewels of numerical analysis
is the QR algorithm for computing the eigenvalues of a matrix.
The core of the QR algorithm is an alternation between a matrix
computation, QR factorization, and an elementary but crucial
nonlinear step, a diagonal shift.  The same can be said of the
steps (2) and (3) that define AAA.}

Step (1) is the feature that distinguishes continuum AAA from
its fixed-grid predecessor.  At each iterative step, the vector
$X$ is determined by the current vector $S$ of support points.
From one step to the next most sample points stay the same, except
in the interval between the two support points where a new support
point $\snew$ has been introduced in step (3); let us call this
interval $(s_j, s_{j+1})$.  At step (1) the three sample points in
$(s_j, s_{j+1})$ are removed and six new sample points are added,
three in $(s_j, \snew)$ and three in $(\snew, s_{j+1})$.

In fact, we have slightly simplified the description of step (1).
Inspection of the code {\tt XS} will show that it takes not one
but two arguments, $X = \hbox{\tt XS}(S,p)$, so that an arbitrary
number $p$ of sample points will be placed between each pair of
support points.  The code {\tt aaax} sets $p = \max\{3, 16-m\}$
at step $m$, so that $p$ begins at $p=14$ for $m=2$ and reduces
linearly with $m$ until it hits $p=3$ for $m\ge 13$.  This is an
engineering adjustment to ensure that the function $f$ is sampled
at dozens of points during the approximation process, rather than
as few as $3$ points.  Such precautions have been been a familiar
feature of adaptive numerical algorithms since the first adaptive
integrators were introduced in the 1960\kern .3pt s.

The other feature of the algorithm to be spelled out is the
convergence criterion.  The iteration stops if the maximum norm
relative error on the sample points falls below a prescribed
value {\tt tol}, set by default to $10^{-13}$.
The actual maximum error over all of $[-1,1]$ will be slightly
larger, and this is checked a posteriori by evaluation on the
finer grid $\hbox{\tt XS}(S,30)$.  (If function evaluations are
expensive, this step can be modified.)  The iteration
also stops if a prescribed maximum degree is reached, which is
set by default to $150$.  In addition the iteration stops if
ten steps in a row have produced ``bad poles,'' that is, one or
more poles in $[-1,1]$, and at least two digits of accuracy have
been achieved.  Poles are computed by means of a generalized
eigenvalue problem by the {\tt prz} code adapted from Chebfun;
this is extremely fast and accurate, adding only about $25\%$
to the overall computation time even though poles are computed
at every AAA step.  In the case of termination due to bad poles,
the approximation $r$ returned by {\tt aaax} is the last one that
was successfully computed without bad poles, so its accuracy will
not be as good as $10^{-13}$.  In practice we often find that a
few digits are lost, but the approximation is always pole-free in the
approximation interval.

\begin{figure}
\begin{center}
\vskip 10pt
\includegraphics[scale=.7]{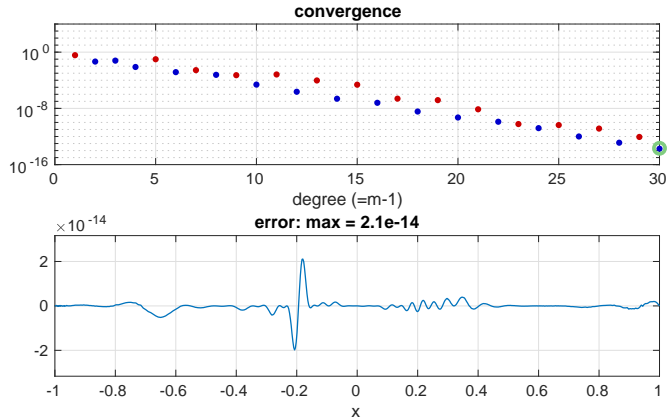}
\end{center}
\caption{\label{fig3} Plots produced by {\tt aaax} in approximating
$f(x) = \exp(-1/x^2)$.
The upper image shows convergence as a function
of degree, with red dots marking steps with ``bad poles'', i.e., poles
in $[-1,1]$.  The red-blue alternation is a result of $f$ being even.
The second image shows the final error curve and
maximal error on a finer grid.  (Computation time excluding plotting:
$20$ ms.)}
\end{figure}

In the default mode of operation, with an invocation as simple
as \verb|r = aaax(f)|, the continuum AAA algorithm produces a
pair of plots showing convergence and the final error curve.
This is illustrated first in Figure~\ref{fig3} for a smooth
function, $f(x) = \exp(-1/x^2)$.
Note that red dots are used to mark
steps of the iteration with bad poles.  The errors are measured
on the sample grid; off the grid, the error in cases with bad
poles would be $\infty$.  The green circle in the upper plot,
corresponding to the error value printed in the title of the lower
plot, is the actual error as computed on the finer plotting grid.

\begin{figure}
\begin{center}
\vskip 10pt
\includegraphics[scale=.7]{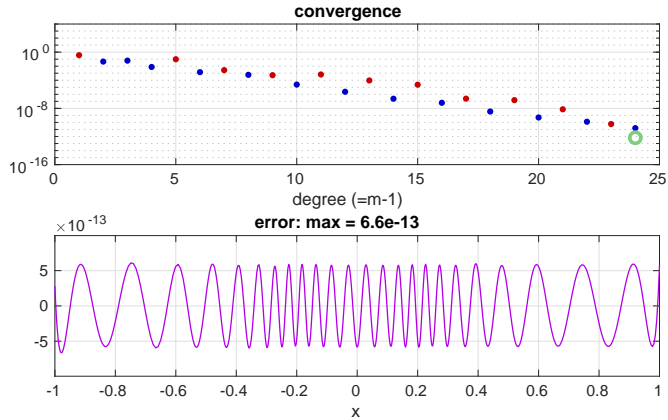}
\end{center}
\caption{\label{fig4}Approximation of the same function $f(x) =
\exp(-1/x^2)$, but now with the function call {\tt aaax(f,24,20)}
to specify maximal degree $24$ and $20$ AAA-Lawson steps.
The result shows the expected equiosillation with an error of
$6.6\times 10^{-13}$ (\/$25$ ms).  The green circle lies below
the final blue dot because AAA-Lawson has improved the error.}
\end{figure}

If a second argument is provided to {\tt aaax}, this is interpreted
as the degree $m-1$ for the approximation (\ref{theapprox}) (more precisely a maximum
degree, since the prescribed tolerance may be achieved earlier).
Typically the result is an approximation whose error is on the
order of a factor of ten or so greater than the minimax error
for that degree.  If a third argument is also provided, this is
interpreted as a number of AAA-Lawson iterative steps to take in an
attempt to improve the approximation to minimax.  The 18 lines of
Lawson code, which can be seen in the Appendix, are adapted from the
algorithm introduced in~\cite{aaalawson}. 
We normally take 20
Lawson steps, and as discussed in~\cite{aaalawson}, convergence
occurs in the majority of cases so long as the error level is
well above machine precision.  Such a computation is illustrated
in Figure~\ref{fig4}, showing the result for degree 24 minimax
approximation of the same function $f(x) = \exp(-1/x^2)$ via
\verb|aaax(@(x) exp(-1./x.^2,24,20))|.
A fourth argument to {\tt aaax} allows one to adjust the convergence
tolerance.

\begin{figure}
\begin{center}
\vskip 10pt
\includegraphics[scale=.7]{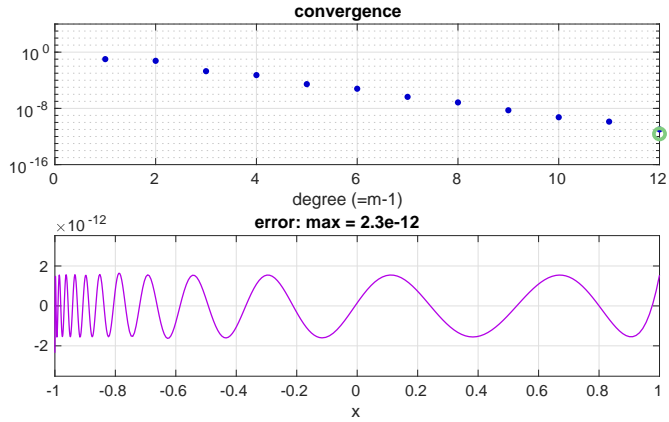}
\end{center}
\caption{\label{fig5}Approximation of\/ $\exp((x-1)/(x+1))$ on 
$[-1,1]$, which is equivalent to the Cody--Meinardus--Varga problem of
approximation of $\exp(s)$ on $(-\infty,0\kern .5pt]$ (\/$20$ ms).  Here
the degree $12$ AAA-Lawson result is computed.}
\end{figure}

Figure~\ref{fig5} shows an example that is equivalent to
a famous rational approximation problem, the degree $n$
rational approximation of $e^{s}$ for $s\in (-\infty,0\kern .3pt]$
first considered by Cody, Meinardus and Varga~\cite{cmv}.  Setting
$s = (x-1)/(x+1)$ transplants this to the problem of degree $n$
approximation of $\exp((x-1)/(x+1))$ for $x\in[-1,1]$.  Executing
\verb|aaax(@(x) exp((x-1)./(x+1))| produces Figure~\ref{fig5} in
about $20$ ms of computing time.  The same approximation can be
computed by Chebfun {\tt minimax}, taking about 100 times as long.
For more on this problem see pp.~214--219 of~\cite{atap}.

The example of Figure~\ref{fig5} provides an illustration of
the importance of the barycentric representation for numerical
stability.  If $r$ is written as a quotient of polynomials $p/q$
for this approximation, then $q$ decreases by a factor
on the order of $10^{11}$ as $x$ moves from $1$ to $-1$, and $p$
decreases by a factor on the order of $10^{23}$.  In floating
point arithmetic, one would need a precision of more than
23 digits to retain any accuracy for $x\approx -1$.
This effect gets rapidly more extreme at higher degrees.

\begin{figure}
\begin{center}
\vskip 10pt
\includegraphics[scale=.7]{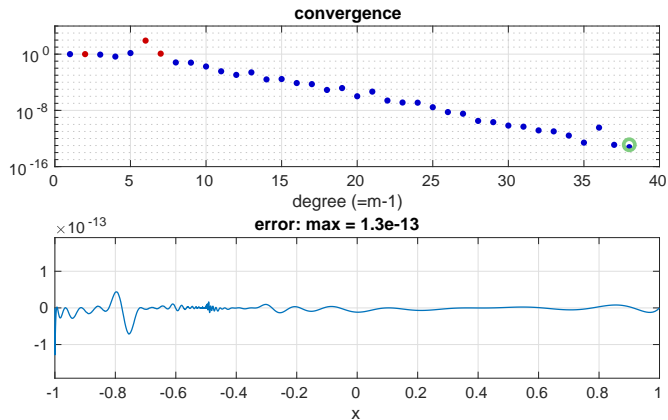}
\end{center}
\caption{\label{fig6}Approximation of the sigmoidal or Fermi--Dirac function
$f(x) = 1/(1+\exp(1000(x+0.5)))$
(\/$20$ ms).}
\end{figure}

Figure~\ref{fig6} considers approximation of the sigmoidal
or Fermi--Dirac function $f(x) = 1/(1+\exp(1000(x+0.5)))$
arising in electronic structure calculations, which makes a rapid
transition from $f(x) \approx 1$ for $x<-0.5$ to $f(x) \approx 0$
for $x> -0.5$.  A degree 38 rational approximation of accuracy
$1.3\times 10^{-13}$ is obtained in 20 ms.\ \ This appears
to offer a big improvement over a recently published algorithm
for this problem~\cite{moussa}.

\begin{figure}
\begin{center}
\vskip 10pt
\includegraphics[scale=.7]{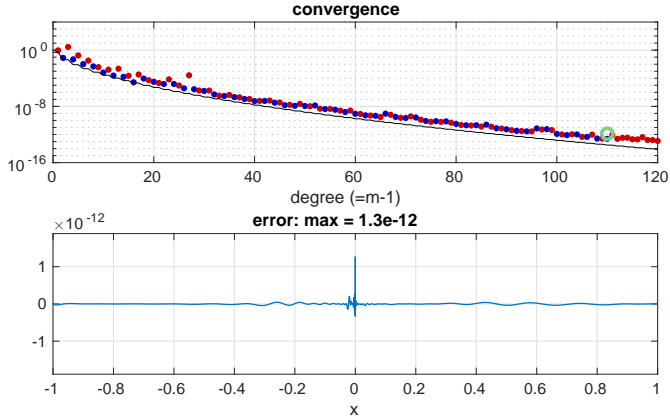}
\end{center}
\caption{\label{fig7}Approximation of the classic function with a
branch point singularity, $f(x) = |x|$, showing root-exponential
convergence as discovered by Newman~{\rm \cite{newman}}.
$12$-digit accuracy is achieved in $600$ ms.\ \ The black
line under the data points shows the error of true minimax
approximation.} \end{figure}

Our final example of this section returns to the problem
shown in Figure~\ref{fig1}, approximation of $f(x) = |x|$.
Figure~\ref{fig7} presents the convergence curve in this case,
showing root-exponential convergence with approximately alternating
blue and red dots (since $f$ is even) until at a level below
$10^{-12}$, all the dots turn red and the iteration terminates.
This successful computation of an approximation accurate to 12
digits in ordinary machine arithmetic happens in less than a
second of computer time.

\begin{figure}
\begin{center}
\vskip 10pt
\includegraphics[scale=.75]{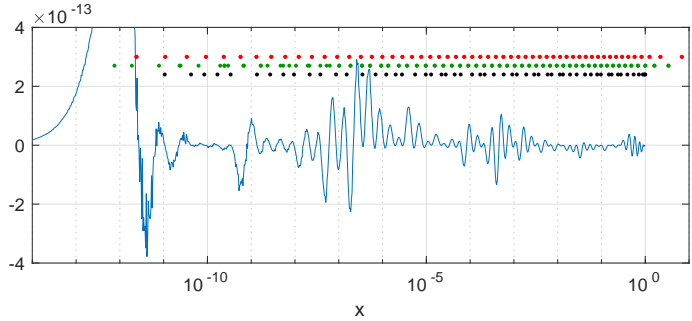}
\end{center}
\caption{\label{fig8}The error $f(x)-r(x)$ in continuum AAA
approximation of $|x|$ on $[-1,1]$, plotted for $x\in (0,1]$ now
on a semilogx scale.  The dots show positions of support points
in $(0,1]$ (black) and absolute values of conjugate pairs of zeros
(green) and poles (red) near the imaginary axis.} \end{figure}

It is known that the power of rational approximation for functions
with singularities derives from exponential clustering of
poles and zeros near these points~\cite{clustering}.  The success
of AAA for such problems depends on the support points being
exponentially clustered too.  Figure~\ref{fig8} illustrates how
the continuum AAA computation of Figure~\ref{fig7} has achieved
the necessary clusterings for all these three sets of points.
The curve in the figure represents the error $f(x)-r(x)$ on the
positive half of the domain, $x\in (0,1]$, and dots are placed on
the same horizontal scale representing support points in $(0,1]$
(black), absolute values of zeros (green), and absolute values
of poles (red).  The poles and zeros lie approximately on the
imaginary axis and in conjugate pairs, so each green and red dot
has multiplicity~$2$.  It is evident that poles and zeros are
(mostly) interlacing and exponentially clustered toward $x=0$,
with the support points exponentially clustered in a similar
manner.  Towards $x=0$ the spacing between dots stretches out, the
``tapering'' effect analyzed in~\cite{clustering}.  All this is
computed blindly by the continuum AAA algorithm, which has no knowledge of
optimality conditions or of the theory of exponential clustering.

Some problems of interest, such as those of
Figures~\ref{fig4} and \ref{fig7}--\ref{fig8}, involve functions
$f(x)$ that are even or odd, and the algorithms could be modified
to exploit this symmetry.  We 
comment on exploitation of symmetry at the end of the next section.

\section{Continuum AAA on the unit circle or disk} For
approximation on the unit circle, our algorithm is mostly the same;
the detailed changes can be seen in the code {\tt aaaz} available in the
Supplementary Materials.\ \ Sample points are now placed between
support points all around the circle with respect to angle, or
equivalently, arc length.  If AAA-Lawson is invoked, a numerical
winding number is calculated since this is of interest for best
and near-best approximations because of Rouch\'e's theorem.
In the simplest case, if a function $f$ analytic in the disk has
a degree $m-1$ rational approximation $r$ that is also analytic in
the disk, and the error curve $(f-r)(|z|=1)$ is nearly a circle of
winding number $\ge 2m-1$, then $\|f-r\|$ must be correspondingly
close to minimal~\cite[Props.~2.1 and 2.2]{nm}.  This is a complex
analogue of the well-known de la Vall\'ee Poussin lower bound in
real approximation on an interval~\cite{atap}.

\begin{figure}
\begin{center}
\vskip 10pt
\includegraphics[scale=.90]{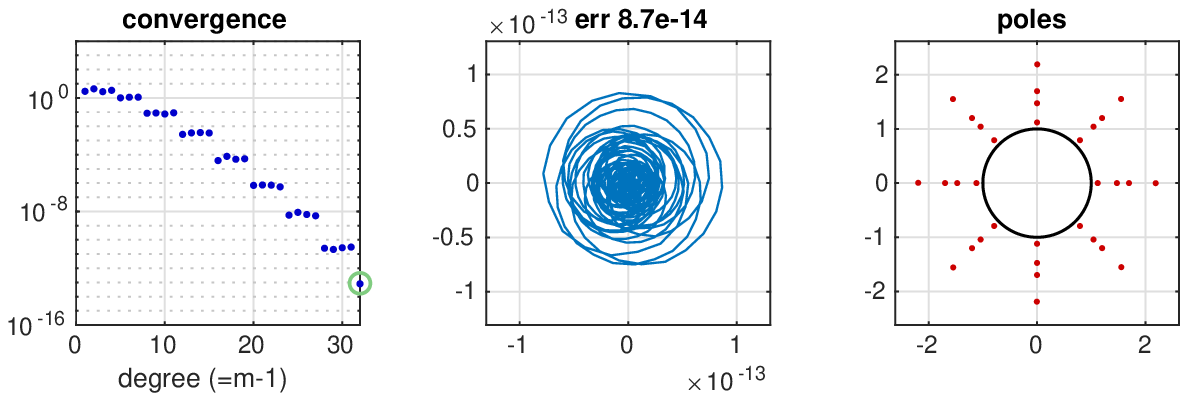}
\vskip 10pt
\includegraphics[scale=.90]{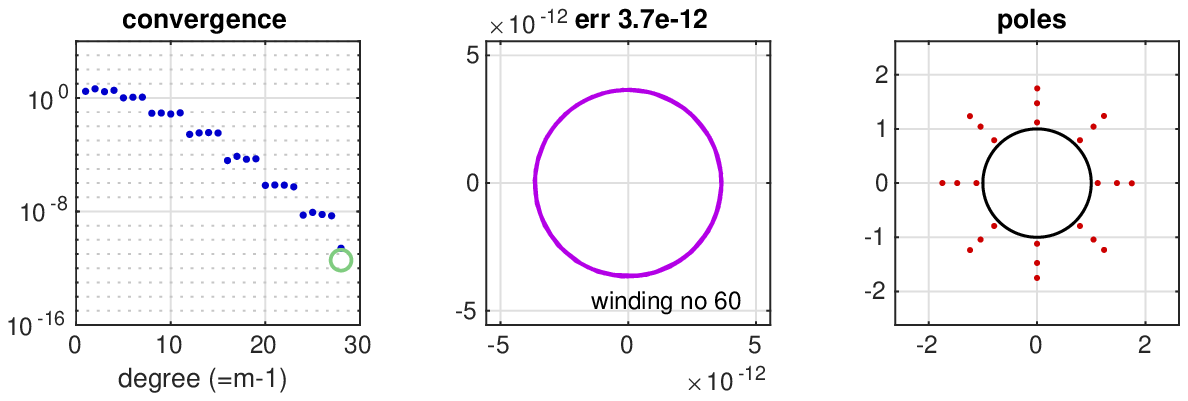}
\end{center}
\caption{\label{fig9}Continuum approximation of $f(z) = \tan(z^4)$  on the unit circle.
Upper row: {\tt aaaz(f)} (\/$30$ ms).
Lower row: {\tt aaaz(f,28,20)} (\/$40$ ms).
In the second case, in addition to the $24$ poles
visible in the plot, there are four more with moduli about $6$.}
\end{figure}

\begin{figure}
\begin{center}
\vskip 10pt
\includegraphics[scale=.90]{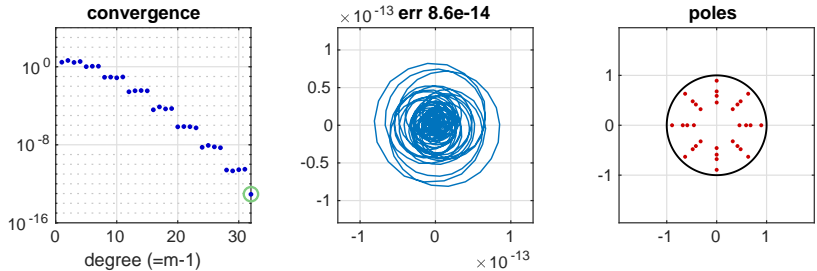}
\end{center}
\caption{\label{fig10}Approximation of $\tan(z^{-4})$ with the {\tt mero}
parameter set to $1$, so that poles of $r$ are permitted in the
interior of the disk (\/$30$ ms).}
\end{figure}

Figure~\ref{fig9} illustrates {\tt aaaz} approximation of
$f(z) = \tan(z^4)$.  This function is analytic in the disk
and meromorphic outside, with eight rays of poles extending
to $\infty$.  Continuum AAA achieves the prescribed accuracy by
approximating the rays each by four poles.  (Note that the errors
in the convergence plot decrease in a staircase pattern because
of the symmetry, a reflection of the Walsh table of best rational
approximants to $f$ breaking into $4\times 4$ blocks of identical
entries~\cite{atap}.)  The innermost eight poles have moduli
matching the value $(\pi/2)^{1/4}$ for $f$ to about 12 digits,
and the next three rings of eight poles match the corresponding
poles of $f$ to about 5, 2, and 1 digits, respectively.  Similarly,
the lower row of the figure shows an approximation to the same
function with degree $28$ and $20$ steps of AAA-Lawson.

Although continuum AAA for the unit circle is like the
algorithm for the unit interval, two new issues arise.  The first
is that the domain now encloses an interior, which raises the
question, should the approximation $r$ be required to be analytic
there?  In other words, if a pole appears with $|z|< 1$, should
it be regarded as a bad pole?  For some problems the answer will
be yes, if an approximation analytic in the disk is sought, and in
other cases it will be no, when we truly wish to approximate just
on the circle.  The code {\tt aaaz} controls this choice with
an input parameter {\tt mero} (``meromorphic''), which by default is set to $0$.
If $\hbox{\tt mero}=0$, then poles in the disk are treated as
bad, following the same logic as described in the last section.
If $\hbox{\tt mero}\ne 0$, then poles in the disk are accepted.
(The possibility of poles exactly on the unit circle is unlikely
enough in floating point arithmetic that we do not worry about it.)

\begin{figure}
\begin{center}
\vskip 10pt
\includegraphics[scale=.90]{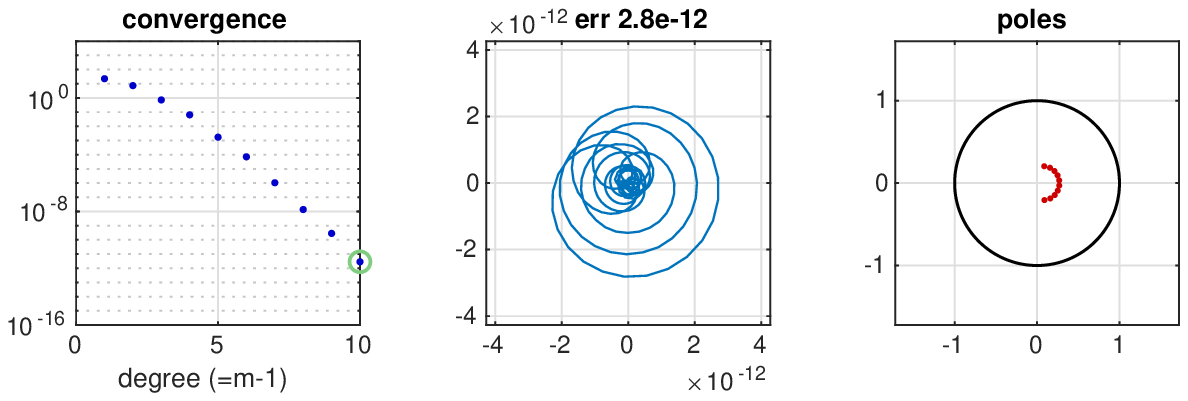}
\vskip 10pt
\includegraphics[scale=.90]{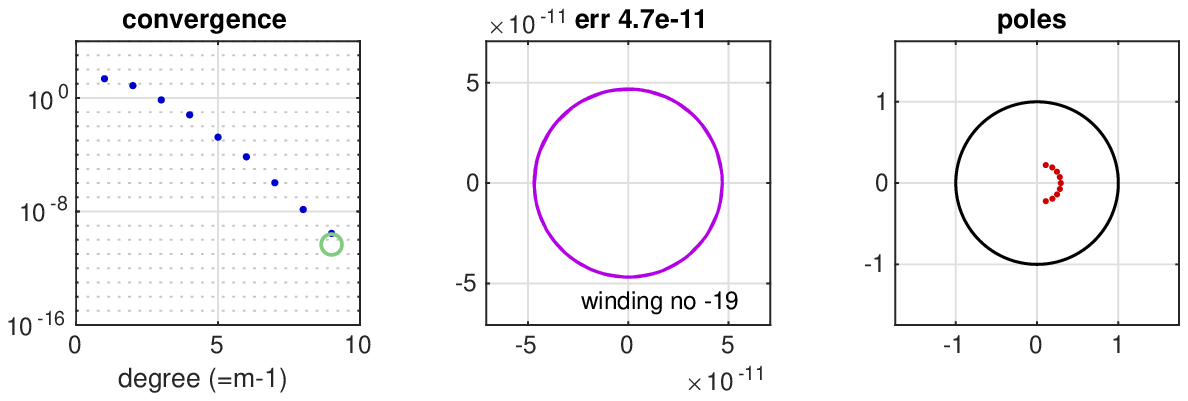}
\end{center}
\caption{\label{fig11}Approximation of
$f(z) = \exp(4/z)$ on the unit circle, with an essential singularity at $z=0$.
Upper row: {\tt aaaz(f)} (\/$30$ ms).
Lower row: {\tt aaaz(f,28,20)} (\/$40$ ms).}
\end{figure}

Suppose, for example, we wish to approximate $f(z) = \tan(z^{-4})$
on the unit circle instead of $\tan(z^4)$.  With the default value
$\hbox{\tt mero}= 0$ this will result in failure (not shown),
and indeed it is obvious that $\|f-r\|$ can be no less than
$\tanh(1) \approx 0.76$ for a function $r$ analytic in the disk
since $\tan(z^{-4})$ has winding number $-4$ on the unit circle
and minimal modulus $\tanh(1)$, whereas an analytic function $r$
must have nonnegative winding number.  (This is Rouch\'e's theorem
again.)  With $\hbox{\tt mero}= 1$, however, the approximation is
straightforward, and the result is shown in Figure~\ref{fig10},
essentially a reflection in the unit circle of the upper row
of Figure~\ref{fig9}.

Figure~\ref{fig11} shows another pair of examples, near-minimax
and minimax, involving the function $f(z) = \exp(4/z)$, which has
an essential singularity at $z=0$.  Although $f$ itself is not
meromorphic in the unit disk, it is still analytic on the unit
circle, and {\tt aaax} with $\hbox{\tt mero}=1$ has no trouble
constructing approximations.

\begin{figure}
\begin{center}
\vskip 10pt
\includegraphics[scale=.90]{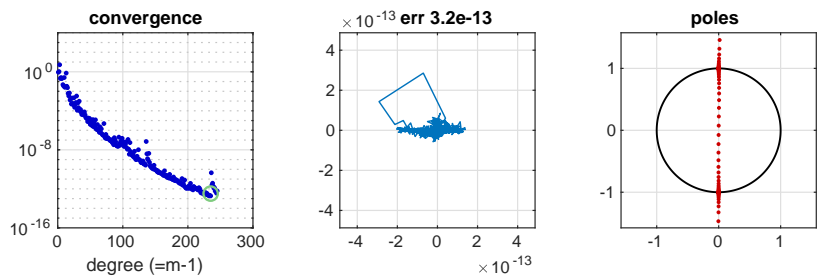}
\end{center}
\caption{\label{fig12}Approximation of\/ $|\Re (z)|$ on the unit
circle with the
{\tt mero} parameter set to $1$, with exponentially clustered
poles on both sides of $z=\pm i$ (\/$9$ s).   Removing
the poles inside or outside of the disk would be the first
step towards a AAA-least squares computation as proposed by
Costa~{\rm\cite{costa2,costa}}.} \end{figure}

Figure~\ref{fig12} shows approximation on the unit circle of $f(z)
= |\Re(z)|$, boundary data that cannot be matched by
any function analytic inside or outside the disk.  In the
fashion familiar in the theories of Wiener--Hopf factorization
and Riemann--Hilbert problems, however, it can still be regarded
as the sum of one function analytic inside the disk plus another
analytic outside (``analytic plus coanalytic'').  As discussed
in~\cite{costa}, 
AAA approximation may provide a valuable tool for such
problems.  The figure shows approximation to 13 digits
by a rational function of degree 235, which has about 60 poles
exponentially clustered on both sides of $z=i$ and $z=-i$.
By separating the two sets of poles and switching to a partial
fractions representation with coefficients determined by linear
least-squares fitting, this can be the basis of very interesting
further computations, including the solution of the Laplace or
Stokes equations~\cite{costa2,xue}.

\begin{figure}
\begin{center}
\vskip 10pt
~~~~~\includegraphics[trim=0 0 0 0, clip, scale=.80]{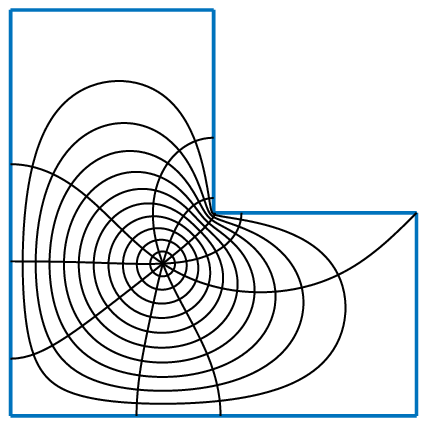}
\vskip 20pt
\includegraphics[scale=.90]{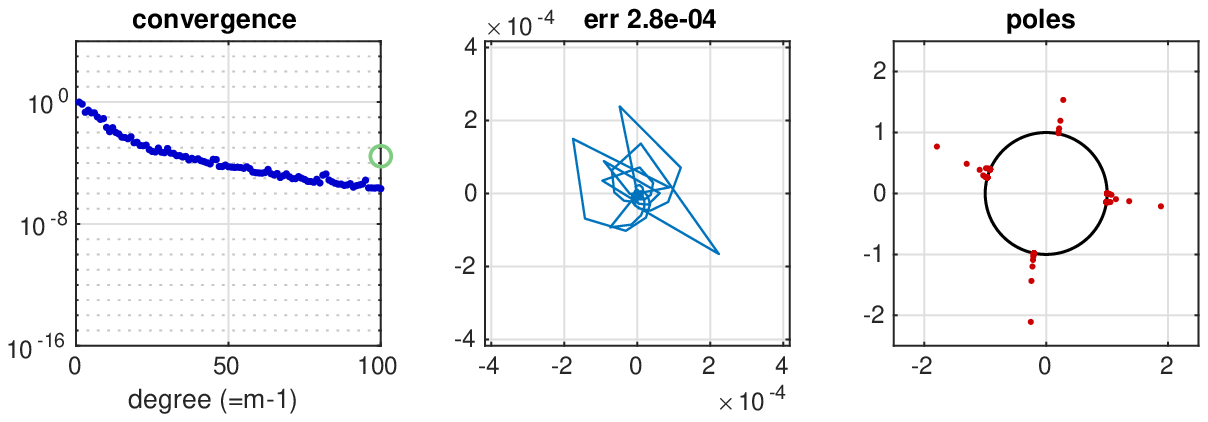}
\end{center}
\caption{\label{fig13}  Conformal mapping example.  The upper plot
shows the conformal image in an L-shaped region of concentric
circles and radial lines in the unit disk, as computed by the
Schwarz--Christoffel Toolbox for MATLAB\/~{\rm \cite{sct}}.
The lower image shows a degree $100$ continuum AAA rational
approximation.  The position of the green circle reflects a loss of accuracy near some of the
vertices.} \end{figure}

Figure~\ref{fig13} returns to functions analytic in the disk with
an example from numerical conformal mapping.  
A conformal map $f$ of the unit disk onto a polygonal region
can be represented as a Schwarz--Christoffel (SC)
integral, and the standard software for computing such maps is the
SC Toolbox for MATLAB~\cite{sct}.  As pointed out in~\cite{cm1},
an effective strategy for representing SC maps numerically is
rational approximation, taking advantage of poles exponentially
clustered near corner singularities.  The figure shows the result
of degree $100$ rational approximation of a map onto an L-shaped
region as generated by essentially this MATLAB code:

{\small
\vskip 10pt
\begin{verbatim}
         v = [0 1i -1+1i -1-1i 1-1i 1];
         f = diskmap(polygon(v)); f = center(f,-.25-.25i); plot(f)
         r = aaaz(f,100);
\end{verbatim}
\vskip 10pt
\par}

\noindent Note that the poles of the approximant cluster
exponentially near six ``prevertices'' along the unit circle.
The accuracy is about 6 digits over most of the domain but
falls to 4 digits very near the vertices (near five of them, to be
precise).  The rational approximation can be used to map points $z$
in the disk to their images in the L shape about ten times faster
than with the SC Toolbox (about 1 ms vs.\ 10 ms per evaluation).
As pointed out in~\cite{cm1}, the ratio of speeds for evaluating
the inverse map $f^{-1}$ typically exceeds $100$.

The other new feature that arises in approximation on the circle
is the matter of real symmetry.  A function $f$ is said to be
{\em real symmetric} (or {\em conjugate symmetric} or {\em Hermitian symmetric})
if $\overline{f(z)} = f(\overline{z})$,
and in such cases one might like to require the approximation $r$
to be real symmetric too.  This is attractive not just cosmetically
but also because exploiting the symmetry can provide a significant
speedup.

Our template code {\tt aaax} does not include an option to
enforce real symmetry, and it breaks symmetry in two ways.
First is simply by rounding errors.
Second and more importantly, even in exact arithmetic, the sequence of support
points will be nonsymmetric in general.  One could modify the code to force symmetry at
all stages, and various authors have done this~\cite{hochman,johns,RF,vre}.
Other types
of symmetry, most notably odd and even symmetry as mentioned at
the end of the last section, could also be optionally enforced.

\section{Continuum AAA on the imaginary axis or right half-plane}
Many applications, especially in control theory and model order reduction,
feature approximations $f(z)\approx r(z)$ on the imaginary axis, typically with a constraint that
the poles should be in the left half-plane: in other words, $r$ should be
analytic in the right half-plane.  These problems
can be regarded as transplants of approximation problems on the unit circle by
a M\"obius transformation.  Specifically, the functions
\begin{equation}
z = M\kern -1pt \left({1+w\over 1-w}\right), \quad w = {z-M\over z+M}
\label{moebius}
\end{equation}
describe a bijection of the unit disk in the $w$-plane and the
right half of the $z$-plane, with $w = -1, 0, 1$ corresponding to
$z = 0, M, \infty$.  As a practical matter, it might be desirable
for a code to include $M$ as a parameter, since in an application
one might happen to know that most of the action of $f(z)$ on
the imaginary axis occurs, say, on the scale $|\Im z| = O(10^6)$
rather than $O(1)$.  Our code arbitrarily fixes $M=1.207$, a choice close to
$1$ but not equal to it so as to lead to fewer surprises since a pole
$z=-M$ goes undetected because it is mapped to $\infty$ by the M\"obius map.

\begin{figure}
\begin{center}
\vskip 10pt
\includegraphics[scale=.75]{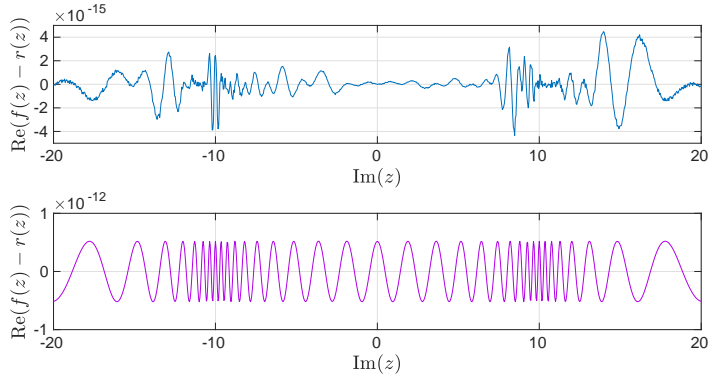}
\end{center}
\caption{\label{fig14}Errors in AAA approximation of
$(\ref{twopts})$ on the imaginary axis to the default tolerance
$10^{-13}$ (\/$25$ ms) and with degree $20$ and $20$ Lawson steps
(\/$70$ ms).  Both $f$ and the error $f-r$ are complex, and just
the real part of the error is plotted.} \end{figure}

Our template code {\tt aaai} is a
shell that calls {\tt aaaz} after the transplantation (\ref{moebius}).
It has the same arguments {\tt f} for the function, {\tt deg} for the 
degree, {\tt nl} for the number of AAA-Lawson iterations, {\tt tol} for
the tolerance, and {\tt mero}
to specify if a meromorphic approximation is permitted (i.e., with poles in the
right half-plane as well as the left).
We have not investigated continuum approximation on the imaginary axis in detail
and give just two examples.  Figure~\ref{fig14} shows near-minimax and minimax
approximations to the function
\begin{equation}
f(z) = {1\over \sqrt{z-a}\kern 1pt \sqrt{z-\overline{a}\kern .5pt}\kern 1pt}, \quad
a = -1 + 10\kern .5pt i.
\label{twopts}
\end{equation}
This function has branch points at $z = -1\pm 10\kern .5pt i$, and
these are reflected in the denser error oscillations 
for $|\Im z| \approx 10$.  The figure shows the real part
of the complex error $f(z)-r(z)$. For the second, minimax
approximation the complex error curve describes a near-circle of winding number $35$ as $z$
traverses the imaginary axis downward from $+\infty\kern .5pt i$
to $-\infty\kern .5pt i$.

\begin{figure}
\begin{center}
\vskip 10pt
\includegraphics[scale=.80]{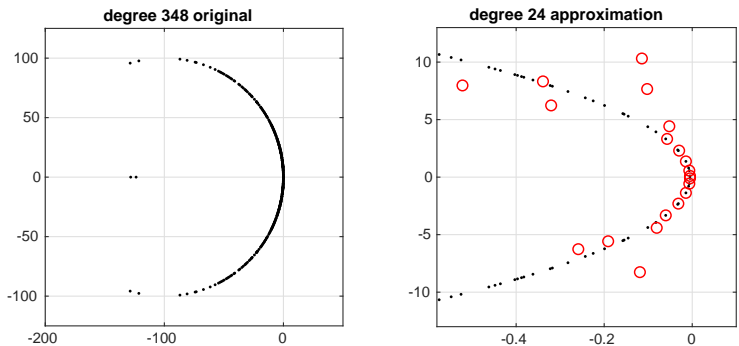}
\end{center}
\caption{\label{fig15}Poles in the complex plane for the clamped
beam example from the NICONET collection~{\rm \cite{vd}}.
On the left, the eigenvalues of the $348\times 348$ matrix $A$,
which are the poles of its resolvent function.  On the right, the
poles of the degree $24$ AAA approximant computed with {\tt aaai}.
The poles near the imaginary axis are closely approximated, giving
about $10^{-4}$ relative accuracy with (if symmetry is exploited)
$72$ function evaluations.} \end{figure}

The second example, shown in Figure~\ref{fig15}, comes from the
``clamped beam'' problem from the NICONET collection of examples
for model order reduction~\cite{vd}, which was also considered in
Figures~6.13--6.14 of the original AAA paper~\cite{aaa}.  The figure
shows the poles of a continuum AAA approximation of degree $24$
with a relative error of about $10^{-4}$.  For this problem,
each function evaluation amounts to a scalarized resolvent,
involving the inverse of a shift of a $348\times 348$ matrix $A$,
and there are $144$ function evaluations.  This number would be
cut in half if real symmetry were exploited, so in effect we are
getting 4-digit accuracy with about $72$ function evaluations as
compared with the number $500$ used in~\cite{aaa}.  The figure
shows most of the 348 eigenvalues of $A$, on the left (a few are
off-scale), and on the right, the 24 poles of the approximation
$r$.  The poles nearest the imaginary axis are captured 
closely,\footnote{More precisely, the rightmost three
conjugate pairs of poles of $r$ match the rightmost, second-rightmost,
and fourth-rightmost conjugate pairs of eigenvalues of $A$ to 
about $6$, $3$, and $3$ digits respectively.  The third-rightmost conjugate
pair of poles of $A$ is not matched by a pole of $r$, and indeed,
this mode appears not to be excited at all in the beam example data.
In particular there is no peak
at the appropriate position $z\approx 0.79\kern .5pt i$ on the imaginary axis in 
Figure~6.13 of~\cite{aaa}.}
and their asymmetric configuration is a reminder that our codes
do not enforce real symmetry.

In this section we have considered approximation on the imaginary
axis, but of course, there are also problems posed on the
real axis.  Analogously, these may come with a requirement of
analyticity in the upper or lower half-plane.  These problems can
be reduced to the former case by a multiplication by $i$ or $-i$.

\section{Other real and complex domains}
The algorithm we have presented extends readily to other domains.
Approximation on $\real$ has just been mentioned.
Approximation on a semiinfinite line like $[\kern .5pt 0,\infty)$
can be treated by a 
M\"obius transformation to $[-1,1]$, just as we used a
transformation to the unit circle for approximation on the
imaginary axis.  Another real domain of interest in applications,
going back to Zolotarev in the 19th century, is a pair of
intervals $[a,b]\cup [c,d\kern .5pt ]$.  Here we would follow
the algorithm essentially as described, with four initial support
points $a,b,c,d$ and new support and sample points introduced in either
$[a,b]$ or $[c,d\kern .5pt]$ at each step.

\begin{figure}
\begin{center}
\vskip 10pt
\includegraphics[scale=.90]{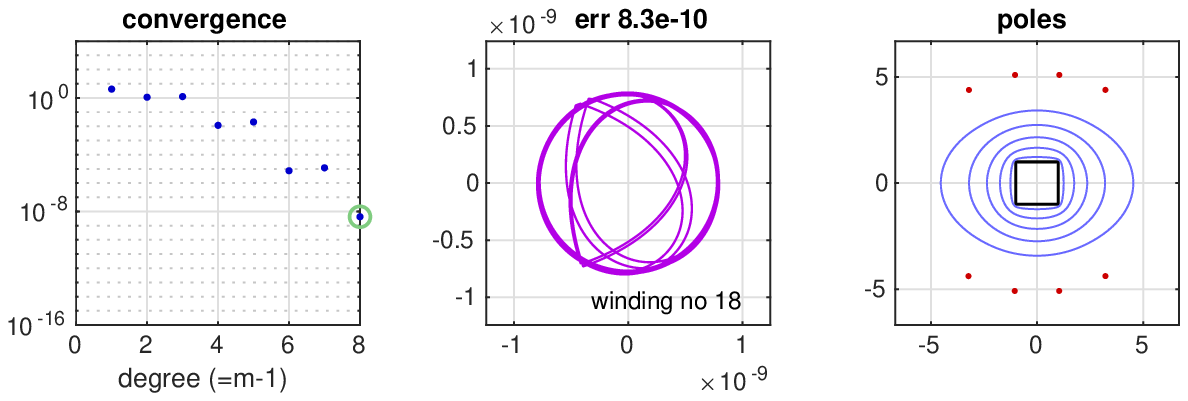}
\vskip 10pt
\includegraphics[scale=.90]{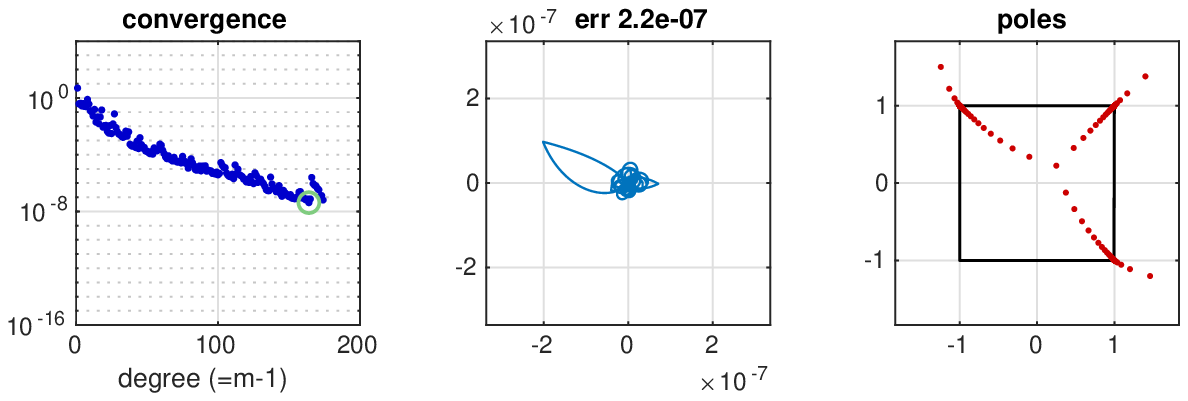}
\end{center}
\caption{\label{fig16} Rational approximations on the unit square.
Upper row: degree $10$ continuum AAA
approximation approximation of $\cos(4z)$ on the square and
its interior with $20$ Lawson steps (\/$100$ ms).  The error curve
is nearly circular apart from four short portions near the corners
of the square, where right angles are preserved since $(f-r)(z)$
is a conformal map.  In the right image, error contours have been
added to shown for $|f(z)-r(z)| = 10^{-8}, 10^{-6}, \dots, 10^0$
(from inside out).  Lower row: continuum AAA approximation of $f(z)
= \min\{\Re(z+1), \Im(z+1)\}$ on the square without interior
(\/$2.7$ ms).  The absence of poles near $z = -1-i$ reveals that
the algorithm has detected that there is no singularity there.}
\end{figure}

On complex domains we follow the pattern for the unit circle
of working on the boundary contour, either with or without a
requirement of analyticity in the interior (i.e., no poles).
No essential change is needed in the algorithm, which can now
distribute sample points according to arc length or
(if the region is starshaped) the angle with respect to a center
point.  For our computed illustrations, we have modified {\tt
aaaz} into a code {\tt aaas} for approximation on the unit square.

Figure~\ref{fig16} shows two examples.  The first consists of
degree $8$ approximation of $\cos(2x)$ on the square with $20$
Lawson steps.  This would take $40$ ms, except that to get a better
image we have run the Lawson iteration on
a grid 10 times finer than usual: 20 SVD calculations involving
matrices of dimension $4509\times 9$.  This reveals an error
curve that is nearly circular apart from four right angles.
In the true minimax error curve, which emerges if one takes 200
rather than $20$ Lawson steps (not shown), the ``double image''
of the figure goes away as symmetry puts two parts of the error
curve in exact superposition, a consequence of $\cos(2z)$ being
an even function.  The corresponding improvement in the error is
from $8.3\times 10^{-10}$ to $8.0 \times 10^{-10}$.

A striking application of rational approximation on complex
domains is to the efficient representation of conformal maps,
as mentioned in Section~3.
In particular, the inverse Schwarz--Christoffel map from a polygonal region to
the unit disk can be represented with great 
efficiency by rational functions~\cite{cm1}.

The second row of Figure~\ref{fig16} shows an example in which the
function to be approximated on the boundary of the square is real,
making it necessary to approximate by a rational function with
poles both inside and outside.  The boundary function chosen,
$f(z) = \min\{\Re(z+1), \Im(z+1)\}$, is zero
on the left and bottom sides.  Thus there is no singularity
at the bottom-left corner, as the zero function is an analytic
continuation to a neighborhood, and the AAA approximation reflects
this in placing no poles near that corner.  The other three
corners have singularities, however, and the poles in the figure
can be interpreted as delineating approximate branch cuts between
the function branches $f_1(z) = 0$ in the lower-left, $f_2(z) =
z+1-i$ near the top, and $f_3(z) = 1+i-iz$ near the right side.
Such approximate branch cut effects are discussed, among
other places, in~\cite{analcont}.
As mentioned in connection with Figure~\ref{fig12}, this kind of
approximation with poles on both sides is the starting point of
the AAA-least squares method introduced by Costa~\cite{costa2,costa}.

\section{Areas for improvement}

The continuum AAA algorithm is remarkably fast and accurate for
most approximation problems.  Like its discrete predecessor, however,
it is sometimes disappointing,
especially in approximation of real functions with
singularities.  Here we discuss four issues with its
behavior, confining our attention to
this most challenging case, approximation by {\tt aaax} of a real function $f(x)$
on $[-1,1]$.

\smallskip

{\em 1. Failure to meet the convergence tolerance.}
When $f$ is smooth, {\tt aaax} typically converges in
milliseconds to the default tolerance $10^{-13}$.  For example,
$f(x) = \tanh(100x)$ is approximated in 20 ms by a degree $30$
rational function with error $1.3\times 10^{-14}$.  With $f(x)
= \tanh(1000x)$, however, red dots appear for
degree $\ge 45$, indicating the presence of bad poles, and the
computation terminates with degree $43$ and error
$1.6\times 10^{-11}$.  This kind of ``red zone termination'' can be
seen in Figure~\ref{fig7}.

\begin{figure}
\begin{center}
\vskip 10pt
\includegraphics[scale=.7]{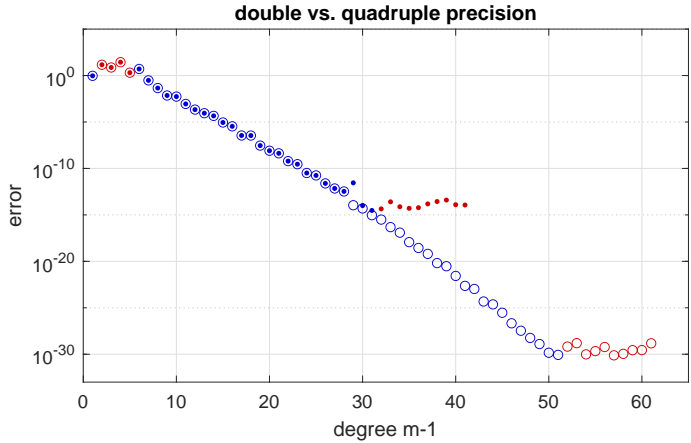}
\end{center}
\caption{\label{fig17}Comparison of double (dots) and quadruple precision
(circles, computed in Julia) for $f(x) = \tanh(100x)$.  Both
iterations are run until ten steps in a row with bad poles have appeared.}
\end{figure}

We believe that in most cases, such failures are related to rounding
errors in floating point arithmetic.  This is confirmed by an
experiment in quadruple precision Julia for $f(x) = \tanh(100x)$,
shown in Figure~\ref{fig17}, where convergence to accuracy
$10^{-29}$ is readily achieved.

Of course, even if a failure would not occur in exact arithmetic,
that does not mean that an algorithm is optimal.  We hope that further
developments will enhance the stability of continuum AAA further,
so that it more reliably gets down to $10^{-13}$ in double precision arithmetic.
\smallskip

{\em 2. Zeroing in on singularities.} Examining the behavior of continuum
AAA for problems with singularities reveals that reasonably uniform
error behavior is often not achieved near these points.  In fact,
trouble near singularities has already appeared in four of our
figures.  In Figures~\ref{fig1} and~\ref{fig7}, both involving $f(x)
= |x|$, there is a disconcerting spike in the final error curve.
In Figure~\ref{fig13} (conformal map to L shape), the same kind of
localized trouble shows up as a green circle two orders of magnitude
above the blue dots, reflecting inaccuracy near prevertex square
root singularities.  It is the semilogx plot of Figure~\ref{fig8}
that reveals the most.  Here we see a much bigger error at the left
side of the plot, off the vertical scale.  The dots in this figure
suggest that AAA would have done better to put one or two support
points closer to the singularity at $x=0$.  We see this in many
experiments involving functions with singularities.  Perhaps a
modification of the algorithm might enhance its ability to zoom
in speedily to difficult points, but we have not yet found the
right trick.  \smallskip

\begin{figure}
\begin{center}
\vskip 10pt
\includegraphics[scale=.75]{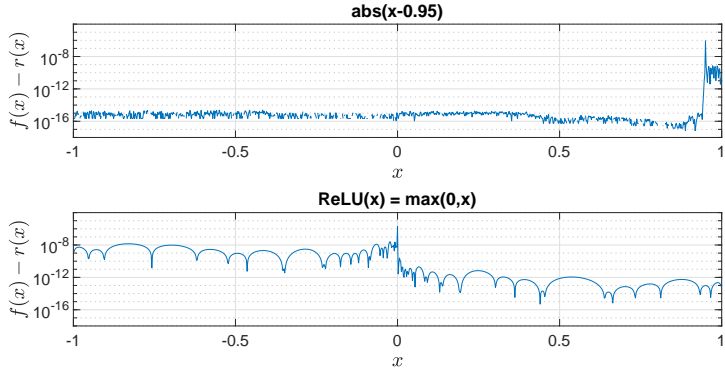}
\end{center}
\caption{\label{fig18}Illustrations of a difficulty that sometimes
arises with the continuum AAA
algorithm: the appearance of imbalances between the two sides of a singularity. The
images show the final errors for $f(x) = |x-0.95|$ and $f(x) = \max(0,x)$.
The errors $7.5\times 10^{-7}$ and $1.5\times 10^{-6}$ are much worse than the error
$1.3\times 10^{-12}$ for the very similar function $f(x) = |x|$.}
\end{figure}

{\em 3. Asymmetries and oscillations.}
Another oddity that arises in certain cases is a pronounced asymmetry
between the errors on the two sides of a singularity, which sometimes
persists on the same side from step to step and sometimes oscillates
from one side to the other (cf.\ Figure~6.1 of~\cite{aaalawson}).
Figure~\ref{fig18} shows two examples, both of which ought to be
straightforward variations of $f(x) = |x|$ but in fact are much
less successful.  In the first plot, for $f(x) = |x-0.95|$, the
final error is ten orders of magnitude higher for $x>0.95$ than for
$x<0.95$.  Clearly something about the algorithm is out of balance.
The second plot shows a similar imbalance in approximation of the
function $f(x) = \max\{ 0 ,x\}$.  The difficulty looks smaller
in magnitude, but it stops convergence just as surely, and it is
particularly embarrassing that our algorithm should do such a poor
job on the famous ``ReLU'' function.

\smallskip

{\em 4. Bad poles.}
Finally there is the perennial question in AAA
approximation---indeed, in rational approximation generally---of
what to do when there are poles in a region
where $r$ should be analytic.  There is an established technology
available for such problems, introduced in the AAA context
in~\cite{hochman} and used by other authors in different settings
(e.g.~\cite{costa,equi}): one can switch from the
barycentric representation to partial fractions, after first
discarding or moving any unwanted poles.  It is also possible
to discard or move poles while retaining a barycentric
representation~\cite{johns1,johns}.

In the early stages of the research that led to this paper, we
assumed that the switch to partial fractions would be an essential
part of our algorithm.  However, in the end we did not find enough
situations where this was advantageous for us to make this a part
of continuum AAA.\ \ One reason to avoid partial fractions is that
they bring a loss of elegance and clarity, since one ends up with
an algorithm that delivers a rational function sometimes in one
representation, sometimes in another.  A more serious drawback is that
there is usually a loss of several digits of precision.  A third
is that, since the poles of a partial fraction representation are
fixed, one loses the possibility of improving a near-minimax fit
to minimax by a Lawson iteration.

Our experiments show that in difficult cases, bad poles often arise
in approximately alternate steps rather than in long sequences
of steps, until one gets down to the ``red zone termination''
related to rounding errors as discussed above.  In other words,
there are usually enough blue dots for the algorithm to proceed to
a successful conclusion.  This need not always be the case, however,
especially for functions with multiple singularities.

Another quite different approach to the bad poles problem is the
``cleanup'' procedure introduced in the original AAA paper~\cite{aaa}, which
has subsequently been refined by the third author and Costa (unpublished).
Here, when a bad pole is encountered, the nearest support point is removed
from $S$ and a new linearized least-squares fit is computed.
Cleanup is not invoked in the algorithm described in the present paper, but perhaps
it should ultimately play a role in the design of robust software.

\section{Summary and discussion}

We have introduced a continuum AAA algorithm for minimax and near-minimax rational approximation
on real and complex sets, with template \hbox{MATLAB} and Julia
programs for $[-1,1]$, the unit circle, and
the imaginary axis.  The algorithm delivers an approximation with no bad poles,
which will usually have relative accuracy at the default tolerance
level of $10^{-13}$.  The speed is remarkable, with approximations typically produced in
milliseconds.    For approximation of degree $m$, the total number
of function evaluations along the way is $\sim 6\kern .7pt m$.

Our template codes lie in the middle of the spectrum from pseudocode to
software.  We hope readers will download them from the Supplementary Materials or the
authors' web sites and enjoy successful explorations.  In the interest of compactness
and readability, however, our codes omit various features that one would expect in
true software.  The omission of an option to impose real symmetry was
discussed in Section 3.  Another example is that  
we have not structured the
codes to avoid recomputation of values $f(z)$ at points $z$ where they have already
been evaluated, which would be important for problems where evaluation of $f$ is expensive.

There does not appear to be much previous work on adaptive selection of
sample points for barycentric
rational approximation, but we note the important recent paper of Pradovera~\cite{pradovera}.
He recommends a strategy in which the next sample point is added at a point
where the barycentric denominator is minimal.

The original, discrete AAA algorithm, whose standard implementation is the {\tt aaa}
code of Chebfun, remains fast and important.  The simplest starting point
in dealing with an approximation problem may
be to simply fix a few thousand points in a set {\tt Z}, evaluate $\hbox{\tt F} = f(\hbox{\tt Z})$, and call
{\tt aaa(F,Z)}.  This is certainly the way to go with exotic sets, such as mixtures
of discrete and continuum components, and it also works for approximating functions that
are meromorphic rather than analytic in the approximation domain, with poles amidst the
sample points.
When one is truly working on a continuum, however,
continuum AAA will be more reliable, both because the adaptive selection of points can
bring more speed and accuracy and also, very importantly, because
of the guarantee that the result will be free of bad poles.

Rational approximation problems are of urgent importance in applied
areas including signal processing and model order reduction, and more recently,
solution of PDE\kern .5pt s.  We are well aware
that this paper is a long
way from such applications, focusing on basic algorithms in the
setting of real and complex approximation theory.

With continuum AAA, is it now feasible to develop a system for 
numerical computation based on rational functions, just as Chebfun
makes use of polynomials and piecewise polynomials?  This is a question
on our minds for ongoing work.

\section*{Acknowledgments}
We are grateful for suggestions from
Stefano Costa, Daan Huybrechs, William Johns,
Davide Pradovera, Mark Reichelt, Olivier S\`ete, Alex Townsend, Michael Tsuk, and Heather
Wilber.
Reichelt and Tsuk are the authors of the rational fitting functionality
in the MathWorks RF Toolbox, based on AAA.

\section{Appendix: MATLAB and Julia code templates}

Here are listings of the main MATLAB template code {\tt aaax} for continuum
AAA approximation on $[-1,1]$ and the three functions {\tt XS},
{\tt prz}, and {\tt reval}, the last two adapted from Chebfun.  These codes
are intended to be read as well as executed, with careful comments in the
mathematical sections but a rougher uncommented style in the plotting sections.

The similar codes {\tt aaaz} and {\tt aaai} for the unit circle and the imaginary axis,
respectively, are available in the Supplementary Materials, together with
Julia equivalents of everything.  The Julia codes offer the option of
computation in quadruple precision.

\indent\vspace*{1pt}

{\small\overfullrule=0pt
\begin{verbatim}
function [r,pol,res,zer,err,S] = aaax(f,deg,nl,tol,plt)
%AAAX   Rational approximation on [-1,1].
%   [r,pol,res,zer,err,S] = aaax(f,deg,nlawson,tol,plt)
%Examples:
%   aaax(@exp); aaax(@exp,5); aaax(@exp,5,20);
%   [r,pol] = aaax(@(x) tanh(20*pi*x)); pol
%   aaax(@abs,150,0,1e-10);

                                               % INITIALIZATIONS
if nargin<2, deg = 150; end                    % arg 2: (max) degree
if nargin<3, nl = 0; end                       % arg 3: Lawson steps, eg 20
if nargin<4, tol = 1e-13; end                  % arg 4: relative tolerance
if nargin<5, plt = 1; end                      % arg 5: plt=0 stops plotting
blue = [0 0 .8]; red = [.8 0 0];
grn = [.5 .8 .5];  MS = 'markersize';
CO = 'color'; LW = 'linewidth'; 
S = [-1 1];                                    % row vector of support pts
f0 = f([S'; XS(S,10)]); err = std(f0);
if (err==0) | (err/abs(mean(f0))<=tol) ...     % constant function or deg 0:
            | (deg==0)                         % return (no plot)
   r = @(x) mean(f0) + 0*x; return
end 
m0 = 2; S0 = S; w0 = [1; -1]; err0 = err;

while true                                     % MAIN AAA LOOP
   m = length(S);
   X = XS(S,max(3,16-m));                      % column vector of sample pts
   C = 1./(X-S);                               % Cauchy matrix
   A = (f(X)-f(S)).*C;                         % Loewner matrix
   [~,~,V] = svd(A,0);                         % SVD
   w = V(:,end);                               % barycentric wt vector
   R = (C*(w.*f(S')))./(C*w);                  % rat approx at sample pts
   err = norm(f(X)-R,inf);                     % abs error of this approx
   pol = prz(S',f(S'),w);                      % poles
   bad = any(imag(pol)==0 & abs(pol)<=1);      % check for bad poles
   co = blue;  
   fmax = norm(f([S';X]),inf);
   if ~bad & err < err0
      m0 = m; S0 = S; w0 = w; err0 = err;      % save latest success
   end
   if plt
      subplot(211), if bad, co = red; end
      semilogy(m-1,err,'.',MS,10,CO,co), hold on
      title convergence, ylim([1e-16 1e4]), grid on
      set(gca,'ytick',10.^(-16:8:0)), xlabel('degree (=m-1)')
   end
   if ~bad & (err/fmax<=tol) | ...             % stop if converged
      (m == deg+1) | ...                       % stop if max degree reached
      ((m-m0>=10) & (err0/fmax<1e-2))          % stop if stagnated
      break
   end 
   [~,j] = max(abs(f(X)-R)); S = [S X(j)];     % next support point
end
m = m0; S = S0; w = w0; err = err0; 
[pol,res,zer] = prz(S',f(S'),w); 
r = @(x) reval(x,S',f(S'),w);                  % function handle for r

if nl > 0                                      % LAWSON ITERATION
   X = XS(S,20); n = length(X);                % switch to finer grid
   wt = ones(m+n,1);                           % Lawson wt vector
   C = 1./(X-S);                               % Cauchy matrix
   R = (C*(w.*f(S')))./(C*w);                  % rat approx at sample pts
   F = [f(X); f(S')]; R = [R; f(S')];          % f and r values
   L = f(X).*C;
   A = [C; eye(m)]; B = [L; diag(f(S))];       % numer and f*denom rows
   for stepno = 1:nl                           % take nl Lawson steps
      [~,~,V] = svd(sqrt(wt).*[A B/fmax],0);
      c = V(:,2*m);
      a = c(1:m); w = -c(m+1:2*m)/fmax;        % numer wts, denom wts
      R = (A*a)./(A*w);
      wt = wt.*abs(F-R);                       % iterative reweighting
      wt = wt/norm(wt,inf);
   end
   r = @(x) reval(x,S',a./w,w);                % function handle for r
end

xx = sort([XS(S,30); S']);                     % COMPUTE ERROR AND PLOT
ee = f(xx)-r(xx); err = norm(ee,inf);
if plt
   semilogy(m-1,max(err,1e-16),'o',LW,2,CO,grn)
   hold off, subplot(212), co = [0 .45 .74];
   if nl > 0, co = [.7 0 .9]; end, plot(xx,ee,LW,.7,CO,co)
   s = sprintf('error: max =%8.1e',err); ylim(1.5*(err+realmin)*[-1 1])
   title(s), grid on, xlabel x
end
end

function XS = XS(S,p)                          % SAMPLE POINTS
S = sort(S); d = (1:p)/(p+1);
XS = S(1:end-1) + d'.*diff(S);
XS = XS(:);
end

function [pol,res,zer] = prz(xj,fj,wj)         % POLES, RESIDUES, ZEROS
m = length(wj);                                %    (from Chebfun)
if any(wj==0)
   ii = find(wj~=0); m = length(ii);
   xj = xj(ii); fj = fj(ii); wj = wj(ii);
end   
B = eye(m+1); B(1,1) = 0;
E = [0 wj.'; ones(m,1) diag(xj)];
pol = eig(E,B); pol = pol(~isinf(pol));
if nargout >= 2
   N = @(t)(1./(t-xj.'))*(fj.*wj);
   Ddiff = @(t) -((1./(t-xj.')).^2)*wj;
   res = N(pol)./Ddiff(pol);
   E = [0 (wj.*fj).'; ones(m,1) diag(xj)];
   zer = eig(E,B); zer = zer(~isinf(zer));
end
end

function r = reval(zz,xj,fj,wj)                % FUNCTION HANDLE FOR R
z = zz(:); C = 1./(z-xj.');                    %    (from Chebfun)    
r = (C*(wj.*fj))./(C*wj);
r(isinf(z)) = sum(wj.*fj)./sum(wj);
ii = find(isnan(r));
for j = 1:length(ii)
   r(ii(j)) = fj(z(ii(j))==xj);
end
r = reshape(r,size(zz));
end
\end{verbatim}
\par}

\end{document}